*Research Article*

# Solutions of Nonlinear Operator Equations by Viscosity Iterative Methods

Mathew Aibinu ,[1] Surendra Thakur,[2] and Sibusiso Moyo[3]

[1]*Institute for Systems Science, Durban University of Technology, Durban 4000, South Africa*
[2]*KZN CoLab, Durban University of Technology, Durban 4000, South Africa*
[3]*Institute for Systems Science & Office of the DVC Research, Innovation & Engagement Milena Court, Durban University of Technology, Durban 4000, South Africa*

Correspondence should be addressed to Mathew Aibinu; moaibinu@yahoo.com

Finding the solutions of nonlinear operator equations has been a subject of research for decades but has recently attracted much attention. This paper studies the convergence of a newly introduced viscosity implicit iterative algorithm to a fixed point of a nonexpansive mapping in Banach spaces. Our technique is indispensable in terms of explicitly clarifying the associated concepts and analysis. The scheme is effective for obtaining the solutions of various nonlinear operator equations as it involves the generalized contraction. The results are applied to obtain a fixed point of $\lambda$-strictly pseudocontractive mappings, solution of $\alpha$-inverse-strongly monotone mappings, and solution of integral equations of Fredholm type.

## 1. Introduction

The Viscosity Approximation Method (VAM) for solving nonlinear operator equations has recently attracted much attention. In 1996, Attouch [1] considered the viscosity solutions of minimization problems. In 2000, Moudafi [2] introduced an explicit viscosity method for nonexpansive mappings. The iterative explicit viscosity sequence $\{x_n\}_{n=1}^{\infty}$ is defined by

$$x_{n+1} = \alpha_n f(x_n) + (1-\alpha_n) T x_n, \quad n \in \mathbb{N}, \quad (1)$$

where $\{\alpha_n\}_{n=1}^{\infty} \subset (0,1)$, $f$ is a contraction on $K$ and the nonexpansive mapping $T: K \longrightarrow K$ is also defined on $K$, which is a nonempty closed convex subset of a Hilbert space $H$. The sequence (1) converges strongly to a fixed point of a nonexpansive mapping $T$ under suitable conditions. Xu et al. [3] proposed the concept of the implicit midpoint rule

$$x_{n+1} = \alpha_n f(x_n) + (1-\alpha_n) T\left(\frac{x_n + x_{n+1}}{2}\right), \quad n \in \mathbb{N}, \quad (2)$$

where $\{\alpha_n\}_{n=1}^{\infty}$, $T$, and $f$ remain as defined in (1). Under certain conditions, they established that the implicit midpoint sequence (2) converges to a fixed point $p$ of $T$ which also solves the variational inequality

$$\langle (I-f)p, x-p \rangle \geq 0, \quad \forall x \in F(T), \quad (3)$$

where $\langle , \rangle$ is the inner product. Aibinu et al. [4] studied the convergence of the sequence (2) in uniformly smooth Banach spaces. Ke and Ma [5] introduced generalized viscosity implicit rules which extend the results of Xu et al. [3]. The generalized viscosity implicit procedures are given by

$$x_{n+1} = \alpha_n f(x_n) + (1-\alpha_n) T(\delta_n x_n + (1-\delta_n) x_{n+1}), \quad n \in \mathbb{N}, \quad (4)$$

$$y_{n+1} = \alpha_n f(y_n) + \beta_n y_n + \gamma_n T(\delta_n y_n + (1-\delta_n) y_{n+1}), \quad n \in \mathbb{N}, \quad (5)$$

where $\{\delta_n\}_{n=1}^{\infty} \subset (0,1)$, $\{\alpha_n\}_{n=1}^{\infty}, \{\beta_n\}_{n=1}^{\infty}, \{\gamma_n\}_{n=1}^{\infty} \subset [0,1]$ with $\alpha_n + \beta_n + \gamma_n = 1$. Replacement of strict contractions in

(5) by the generalized contractions and extension to uniformly smooth Banach spaces was considered by Yan et al. [6]. Under certain conditions imposed on the parameters involved, the sequence $\{x_n\}_{n=1}^{\infty}$ converges strongly to a fixed point $p$ of the nonexpansive mapping $T$, which is also the unique solution of the variational inequality

$$\langle (I-f)p, J(x-p) \rangle \geq 0, \quad \text{for all } x \in F(T), \quad (6)$$

where $J$ is a normalized duality mapping and $\langle , \rangle$ is the duality pairing. Aibinu and Kim [7] used an analytical method to compare the rate of convergence of the sequences (4) and (5). Due to the important roles which nonlinear operator equations play in modeling many phenomena in scientific fields, research on the solution and application of nonlinear operator equations are ever green (see e.g., [8–12]).

Inspired by the previous works in this direction, this paper studies an implicit iterative sequence that involves the generalized contraction and which is effective for obtaining the solutions of various nonlinear operator equations. Precisely, for a nonempty closed convex subset $K$ of a uniformly smooth Banach space $E$ and real sequences $\{\{\alpha_n^i\}_{n=1}^{\infty}\}_{i=1}^{3} \subset [0,1]$ and $\{\delta_n\}_{n=1}^{\infty} \subset (0,1)$ such that $\sum_{i=1}^{3} \alpha_n^i = 1$, an implicit iterative scheme is defined from an arbitrary $x_1 \in K$ by

$$x_{n+1} = \alpha_n^1 f(x_n) + \alpha_n^2 x_n + \alpha_n^3 T((1-\delta_n)f(x_n) + \delta_n x_{n+1}), \quad (7)$$

where $T : K \longrightarrow K$ is a nonexpansive mapping and $f$ is a generalized contraction mapping. The sequence is shown to converge strongly to a fixed point of a nonexpansive mapping in Banach spaces. The adopted technique is indispensable in terms of explicitly clarifying the associated concepts and analysis. The results are applied to obtain a fixed point of $\lambda$-strictly pseudocontractive mapping, solution of $\alpha$-inverse-strongly monotone mapping, and solution of the integral equation of Fredholm type. An example of real sequences which satisfy the stated conditions of our iteration is given.

## 2. Preliminaries

*Definition 1.* Let $E$ be a real Banach space with dual $E^*$ and denote the norm on $E$ by $\|.\|$. The normalized duality mapping $J : E \longrightarrow 2^{E^*}$ is defined as

$$J(x) = \{f \in E^* : \langle x, f \rangle = \|x\| \|f\|, \|x\| = \|f\|\}, \quad (8)$$

where $\langle ., . \rangle$ is the duality pairing between $E$ and $E^*$.

*Definition 2.* Let $B_E \coloneqq \{x \in E : \|x\| = 1\}$. $E$ is said to be smooth (or Gâteaux differentiable) if the limit

$$\lim_{t \to 0^+} \frac{\|x + ty\| - \|x\|}{t} \quad (9)$$

exists for each $x, y \in B_E$. $E$ is said to have uniformly Gâteaux differentiable norm if for each $y \in B_E$, the limit is attained uniformly for $x \in B_E$ and uniformly smooth if it is smooth and the limit is attained uniformly for each $x, y \in B_E$. Also, $E$ is said to be uniformly smooth if

$$\lim_{\tau \to 0} \frac{\rho_E(\tau)}{\tau} = 0. \quad (10)$$

It is well known that $\rho_E$ is nondecreasing. $E$ is said to be $q$-uniformly smooth if there exist a constant $c > 0$ and a real number $q > 1$ such that $\rho_E(\tau) \leq c\tau^q$. $L_p$, $l_p$, and $W_p^m (1 < p < \infty)$ are typical examples of such spaces, where

$$L_p, l_p, W_p^m \text{ is } \begin{cases} 2 - \text{uniformly smooth} & \text{if } 2 \leq p < \infty, \\ p - \text{uniformly smooth} & \text{if } 1 < p < 2. \end{cases} \quad (11)$$

Recall that if $E$ is smooth, then $J$ is single valued and onto if $E$ is reflexive. Furthermore, the normalized duality mapping $J$ is uniformly continuous on bounded subsets of $E$ from the strong topology of $E$ to the weak-star topology of $E^*$ if $E$ is a Banach space with a uniformly Gâteaux differentiable norm.

*Definition 3.* Let $(E, d)$ be a metric space and $K$ a subset of $E$ with $f : K \longrightarrow K$ a mapping defined on $K$.

(i) $f$ is said to be Lipschitzian if there exists a constant $c > 0$, such that for all $x, y \in K$

$$d(fx, fy) \leq cd(x, y). \quad (12)$$

$f$ is called nonexpansive if $c = 1$, and it is a contraction if $c \in [0, 1)$. A contraction mapping $f$ will be referred to as *c-contraction* mapping

(ii) $f$ is said to be a *Meir-Keeler contraction* if for each $\varepsilon > 0$ there exists $\delta = \delta(\varepsilon) > 0$ such that for each $x, y \in K$, with $\varepsilon \leq d(x, y) < \varepsilon + \delta$, we have $d(f(x), f(y)) < \varepsilon$

(iii) Let $\mathbb{N}$ be the set of all positive integers and $\mathbb{R}^+$ the set of all positive real numbers. A mapping $\psi : \mathbb{R}^+ \longrightarrow \mathbb{R}^+$ is said to be an *L-function* if $\psi(0) = 0$, $\psi(t) > 0$ for all $t > 0$ and for every $s > 0$, there exists $u > s$ such that $\psi(t) \leq s$ for each $t \in [s, u]$

(iv) $f : E \longrightarrow E$ is called a $(\psi, L)$-*contraction* if $\psi : \mathbb{R}^+ \longrightarrow \mathbb{R}^+$ is an *L-function* and $d(f(x), f(y)) < \psi(d(x, y))$, for all $x, y \in E, x \neq y$

Throughout this paper, the generalized contraction mappings will refer to Meir-Keeler or $(\psi, L)$-contraction contractions. It is assumed that the *L*-function from the definition of $(\psi, L)$-contraction is continuous, strictly increasing, and $\lim_{t \to \infty} \phi(t) = \infty$, where $\phi(t) = t - \psi(t)$ for all $t \in \mathbb{R}^+$. Whenever there is no confusion, $\phi(t)$ and $\psi(t)$ will be written as $\phi t$ and $\psi t$, respectively.

We have the following interesting results about the Meir-Keeler contraction.

**Proposition 4 Meir and Keeler** [13]. *Let $(E, d)$ be a complete metric space and let $f$ be a Meir-Keeler contraction on $E$. Then, $f$ has a unique fixed point in $E$.*

**Proposition 5 Suzuki** [14]. *Let $E$ be a Banach space, $K$ a convex subset of $E$ and $f : K \longrightarrow K$ a Meir-Keeler contraction. Then, $\forall \varepsilon > 0$, there exists $c \in (0, 1)$ such that*

$$\|f(u) - f(v)\| \leq c\|u - v\|, \tag{13}$$

*for all $u, v \in K$ with $\|u - v\| \geq \varepsilon$.*

**Proposition 6 Lim** [15]. *Let $(E, d)$ be a metric space and $f : E \longrightarrow E$ be a mapping. The following assertions are equivalent:*

(i) *$f$ is a Meir-Keeler type mapping*

(ii) *there exists an L-function $\psi : \mathbb{R}^+ \longrightarrow \mathbb{R}^+$ such that $f$ is a $(\psi, L)$-contraction*

**Proposition 7 Lim** [15]. *Let $K$ be a nonempty convex subset of a Banach space $E$, $T : K \longrightarrow K$ a nonexpansive mapping and $f : K \longrightarrow K$ a Meir-Keeler contraction. Then, $Tf$ and $fT$ are Meir-Keeler contractions.*

The following lemmas are needed in the sequel.

**Lemma 8 Suzuki** [16]. *Let $\{u_n\}_{n=1}^{\infty}$ and $\{v_n\}_{n=1}^{\infty}$ be bounded sequences in a Banach space $E$ and $\{t_n\}_{n=1}^{\infty}$ be a sequence in $[0, 1]$ with $0 < \liminf_{n \to \infty} t_n \leq \limsup_{n \to \infty} t_n < 1$. Suppose that $u_{n+1} = (1 - t_n)u_n + t_n v_n$ for all $n \geq 0$ and $\limsup_{n \to \infty}(\|u_{n+1} - u_n\| - \|v_{n+1} - v_n\|) \leq 0$. Then, $\lim_{n \to \infty} \|u_n - v_n\| = 0$.*

**Lemma 9 Sunthrayuth and Kumam** [17]. *Let $K$ be a nonempty closed and convex subset of a uniformly smooth Banach space $E$. Let $T : K \longrightarrow K$ be a nonexpansive mapping such that $F(T) \neq \emptyset$ and $f : K \longrightarrow K$ be a generalized contraction mapping. Assume that $\{x_t\}$ defined by $x_t = tf(x_t) + (1 - t)Tx_t$ for $t \in (0, 1)$ converges strongly to $p \in F(T)$ as $t \longrightarrow 0$. Suppose that $\{x_n\}$ is a bounded sequence such that $\|x_n - Tx_n\| \longrightarrow 0$ as $n \longrightarrow \infty$. Then,*

$$\limsup_{n \to \infty} \langle f(p) - p, J(x_n - p) \rangle \leq 0. \tag{14}$$

**Lemma 10 Sunthrayuth and Kumam** [17]. *Let $K$ be a nonempty closed and convex subset of a uniformly smooth Banach space $E$. Let $T : K \longrightarrow K$ be a nonexpansive mapping such that $F(T) \neq \emptyset$ and $f : K \longrightarrow K$ be a generalized contraction mapping. Then, $\{x_t\}$ defined by $x_t = tf(x_t) + (1 - t)Tx_t$ for $t \in (0, 1)$ converges strongly to $p \in F(T)$, which solves the following variational inequality:*

$$\langle f(p) - p, J(z - p) \rangle \leq 0, \quad \forall z \in F(T). \tag{15}$$

**Lemma 11 Xu** [18]. *Let $\{a_n\}$ be a sequence of nonnegative real numbers satisfying the following relations:*

$$a_{n+1} \leq (1 - \alpha_n)a_n + \alpha_n \sigma_n + \gamma_n, \quad n \in \mathbb{N}, \tag{16}$$

*where*

(i) *$\{\alpha\}_n \subset (0, 1)$, $\sum_{n=1}^{\infty} \alpha_n = \infty$*

(ii) *$\limsup\{\sigma\}_n \leq 0$*

(iii) *$\gamma_n \geq 0$, $\sum_{n=1}^{\infty} \gamma_n < \infty$*

*Then, $\lim_{n \to \infty} a_n = 0$.*

## 3. Main Results

*Assumption 12.* Let $K$ be a nonempty closed convex subset of a uniformly smooth Banach space $E$ and $f : K \longrightarrow K$ a generalized contraction mapping. Let $T$ be a nonexpansive self-mapping defined on $K$ with $F(T) \neq \emptyset$. The real sequences $\{\{\alpha_n^i\}_{n=1}^{\infty}\}_{i=1}^{3} \subset [0, 1]$ and $\{\delta_n\}_{n=1}^{\infty} \subset (0, 1)$ are assumed to satisfy the following conditions:

(i) $\sum_{i=1}^{3} \alpha_n^i = 1$

(ii) $\lim_{n \to \infty}(1 - \alpha_n^3 \delta_n - \alpha_n^2) = 0$, $\sum_{n=1}^{\infty}(1 - \alpha_n^3 \delta_n - \alpha_n^2) = \infty$

(iii) $0 < \liminf_{n \to \infty} \alpha_n^2 \leq \limsup_{n \to \infty} \alpha_n^2 < 1$

(iv) $\lim_{n \to \infty} \alpha_n^3 = 0$, $\sum_{n=1}^{\infty} \alpha_n^3 (1 - \delta_n) < \infty$

(v) $0 < \varepsilon \leq \delta_n \leq \delta_{n+1} \leq \delta < 1$, for all $n \in \mathbb{N}$

Under the conditions (i)-(v) of Assumption 12 stated above, this study establishes the convergence of the iterative scheme (7).

Firstly, it is shown that for all $\omega \in K$, the mapping defined by

$$u \mapsto T_\omega(u) \coloneqq \alpha_n^1 f(\omega) + \alpha_n^2 \omega + \alpha_n^3 T((1 - \delta_n)f(\omega) + \delta_n u), \tag{17}$$

for all $u \in K$, where $\{\{\alpha_n\}_{n=1}^{\infty}\}_{i=1}^{3} \subset [0, 1]$, $\{\delta_n\}_{n=1}^{\infty} \subset (0, 1)$, is a contraction with $\delta \in (0, 1)$ a contractive constant.

Indeed, for all $u, v \in K$,

$$\begin{aligned}
\|T_\omega(u) - T_\omega(v)\| &= \alpha_n^3 \|T((1 - \delta_n)f(\omega) + \delta_n u) \\
&\quad - T((1 - \delta_n)f(\omega) + \delta_n v)\| \\
&\leq \alpha_n^3 \|(1 - \delta_n)f(\omega) + \delta_n u \\
&\quad - (1 - \delta_n)f(\omega) - \delta_n v\| \\
&\leq \alpha_n^3 \delta_n \|u - v\| \leq \delta_n \|u - v\| \\
&\leq \delta \|u - v\|.
\end{aligned} \tag{18}$$

Therefore, $T_\omega$ is a contraction. By Banach's contraction mapping principle, $T_\omega$ has a fixed point.

The proof of the following lemmas which are useful in establishing the main result are given.

**Lemma 13.** *Let $K$ be a nonempty closed convex subset of a uniformly smooth Banach space $E$ and $f : K \longrightarrow K$ a generalized contraction mapping. Let $T$ be a nonexpansive self-mapping defined on $K$ with $F(T) \neq \emptyset$. For an arbitrary $x_1 \in K$, define the iterative sequence $\{x_n\}_{n=1}^{\infty}$ by (7). Then, the sequence $\{x_n\}_{n=1}^{\infty}$ is bounded under the conditions (i)-(v) of Assumption 12.*

*Proof.* It is needed to show that the sequence $\{x_n\}_{n=1}^{\infty}$ is bounded. For $p \in F(T)$,

$$\begin{aligned}
\|x_{n+1} - p\| &= \|\alpha_n^1 f(x_n) + \alpha_n^2 x_n + \alpha_n^3 T((1 - \delta_n) f(x_n) \\
&\quad + \delta_n x_{n+1}) - p\| \leq \alpha_n^1 \|f(x_n) - p\| + \alpha_n^2 \|x_n - p\| \\
&\quad + \alpha_n^3 \|T((1 - \delta_n) f(x_n) + \delta_n x_{n+1}) - p\| \\
&\leq \alpha_n^1 \|f(x_n) - f(p)\| + \alpha_n^1 \|f(p) - p\| + \alpha_n^2 \|x_n - p\| \\
&\quad + \alpha_n^3 \|(1 - \delta_n) f(x_n) + \delta_n x_{n+1} - p\| \\
&= \alpha_n^1 \|f(x_n) - f(p)\| + \alpha_n^1 \|f(p) - p\| + \alpha_n^2 \|x_n - p\| \\
&\quad + \alpha_n^3 \|(1 - \delta_n)(f(x_n) - p) + \delta_n (x_{n+1} - p)\| \\
&\leq \alpha_n^1 \|f(x_n) - f(p)\| + \alpha_n^1 \|f(p) - p\| + \alpha_n^2 \|x_n - p\| \\
&\quad + \alpha_n^3 (1 - \delta_n) \|f(x_n) - f(p)\| + \alpha_n^3 (1 - \delta_n) \\
&\quad \cdot \|f(p) - p\| + \alpha_n^3 \delta_n \|x_{n+1} - p\| \\
&\leq \alpha_n^1 \psi \|x_n - p\| + \alpha_n^1 \|f(p) - p\| + \alpha_n^2 \|x_n - p\| \\
&\quad + \alpha_n^3 (1 - \delta_n) \psi \|x_n - p\| + \alpha_n^3 (1 - \delta_n) \|f(p) - p\| \\
&\quad + \alpha_n^3 \delta_n \|x_{n+1} - p\| = (\alpha_n^1 \psi + \alpha_n^2 + \alpha_n^3 (1 - \delta_n) \psi) \\
&\quad \cdot \|x_n - p\| + (\alpha_n^1 + \alpha_n^3 (1 - \delta_n)) \|f(p) - p\| \\
&\quad + \alpha_n^3 \delta_n \|x_{n+1} - p\| = ((\alpha_n^1 + \alpha_n^3) \psi + \alpha_n^2 - \alpha_n^3 \delta_n \psi) \\
&\quad \cdot \|x_n - p\| + ((\alpha_n^1 + \alpha_n^3) - \alpha_n^3 \delta_n) \|f(p) - p\| \\
&\quad + \alpha_n^3 \delta_n \|x_{n+1} - p\| = ((1 - \alpha_n^2) \psi + \alpha_n^2 - \alpha_n^3 \delta_n \psi) \\
&\quad \cdot \|x_n - p\| + (1 - \alpha_n^2 - \alpha_n^3 \delta_n) \|f(p) - p\| \\
&\quad + \alpha_n^3 \delta_n \|x_{n+1} - p\| = (\psi + \alpha_n^2 (1 - \psi) - \alpha_n^3 \delta_n \psi) \\
&\quad \cdot \|x_n - p\| + (1 - \alpha_n^2 - \alpha_n^3 \delta_n) \|f(p) - p\| \\
&\quad + \alpha_n^3 \delta_n \|x_{n+1} - p\|.
\end{aligned} \tag{19}$$

Therefore,

$$\begin{aligned}
\|x_{n+1} - p\| &\leq \frac{\psi + \alpha_n^2 (1 - \psi) - \alpha_n^3 \delta_n \psi}{1 - \alpha_n^3 \delta_n} \|x_n - p\| \\
&\quad + \frac{1 - \alpha_n^2 - \alpha_n^3 \delta_n}{1 - \alpha_n^3 \delta_n} \|f(p) - p\| \\
&\quad + \frac{1 - \alpha_n^2 - \alpha_n^3 \delta_n}{1 - \alpha_n^3 \delta_n} \|f(p) - p\| \\
&= \left(1 - \frac{(1 - \alpha_n^2 - \alpha_n^3 \delta_n) \phi}{1 - \alpha_n^3 \delta_n}\right) \|x_n - p\| \\
&\quad + \frac{(1 - \alpha_n^2 - \alpha_n^3 \delta_n) \phi}{1 - \alpha_n^3 \delta_n} \phi^{-1} \|f(p) - p\| \\
&\leq \max \{\|x_n - p\|, \phi^{-1} \|f(p) - p\|\}.
\end{aligned} \tag{20}$$

Then by induction, we have

$$\|x_{n+1} - p\| \leq \max \{\|x_1 - p\|, \phi^{-1} \|f(p) - p\|\}. \tag{21}$$

For $p \in F(T)$,

$$\begin{aligned}
\|f(x_n)\| &\leq \|f(x_n) - f(p)\| + \|f(p)\| \\
&\leq \psi \|x_n - p\| + \|f(p)\| \\
&\leq \max \{\psi \|x_1 - p\|, \psi \phi^{-1} \|f(p) - p\|\} \\
&\quad + \|f(p)\| \text{(by induction)}.
\end{aligned} \tag{22}$$

So, $\{x_n\}_{n=1}^{\infty}$ is bounded. Also,

$$\begin{aligned}
\|T((1 &- \delta_n) f(x_n) + \delta_n x_{n+1})\| \\
&= \|T((1 - \delta_n) f(x_n) + \delta_n x_{n+1}) - p + p\| \\
&\leq \|T((1 - \delta_n) f(x_n) + \delta_n x_{n+1}) - Tp\| + \|p\| \\
&\leq \|(1 - \delta_n) f(x_n) + \delta_n x_{n+1} - p\| + \|p\| \\
&\leq (1 - \delta_n) \|f(x_n) - p\| + \delta_n \|x_{n+1} - p\| + \|p\| \\
&\leq (1 - \delta_n) \|f(x_n) - f(p)\| + (1 - \delta_n) \|f(p) - p\| \\
&\quad + \delta_n \|x_{n+1} - p\| + \|p\| \leq (1 - \delta_n) \psi \|x_n - p\| \\
&\quad + \delta_n \|x_{n+1} - p\| + (1 - \delta_n) \|f(p) - p\| + \|p\| \\
&\leq (1 - \varepsilon) \psi \|x_n - p\| + \delta \|x_{n+1} - p\| \\
&\quad + (1 - \varepsilon) \|f(p) - p\| + \|p\|.
\end{aligned} \tag{23}$$

Therefore,

$$\begin{aligned}
\|T((1 &- \delta_n) f(x_n) + \delta_n x_{n+1})\| \\
&\leq (1 + \delta - \varepsilon \psi) \max \{\|x_n - p\|, \phi^{-1} \|f(p) - p\|\} \\
&\quad + (1 - \varepsilon) \|f(p) - p\| + \|p\| \text{(by induction)}.
\end{aligned} \tag{24}$$

Hence, $\{T((1 - \delta_n) f(x_n) + \delta_n x_{n+1})\}_{n=1}^{\infty}$ is bounded.

**Lemma 14.** *Let $K$ be a nonempty closed convex subset of a uniformly smooth Banach space $E$ and $f : K \longrightarrow K$ a generalized contraction mapping. Let $T$ be a nonexpansive self-mapping defined on $K$ with $F(T) \neq \emptyset$. Suppose $\{\delta_n\}_{n=1}^{\infty}$ is a real sequences in $(0,1)$ and $\{x_n\}_{n=1}^{\infty} \subset K$. Set $y_n = (1 - \delta_n) f(x_n) + \delta_n x_{n+1}$, then*

$$\begin{aligned}
\|Ty_{n+1} - Ty_n\| &\leq (1 - \delta_{n+1}) \psi \|x_{n+1} - x_n\| + (\delta_{n+1} - \delta_n) \\
&\quad \cdot \|x_{n+1} - f(x_n)\| + \delta_{n+1} \|x_{n+2} - x_{n+1}\|.
\end{aligned} \tag{25}$$

*Proof.*

$$\begin{aligned}
&\|Ty_{n+1} - Ty_n\| \\
&= \|T((1-\delta_{n+1})f(x_{n+1}) + \delta_{n+1}x_{n+2}) - T((1-\delta_n)f(x_n) \\
&\quad + \delta_n x_{n+1})\| \leq \|(1-\delta_{n+1})f(x_{n+1}) + \delta_{n+1}x_{n+2} \\
&\quad - (1-\delta_n)f(x_n) - \delta_n x_{n+1}\| = \|(1-\delta_{n+1})f(x_{n+1}) \\
&\quad - (1-\delta_{n+1})f(x_n) + (1-\delta_{n+1})f(x_n) - (1-\delta_n)f(x_n) \\
&\quad + \delta_{n+1}x_{n+2} - \delta_{n+1}x_{n+1} + \delta_{n+1}x_{n+1} - \delta_n x_{n+1}\| \\
&= \|(1-\delta_{n+1})(f(x_{n+1}) - f(x_n)) - (\delta_{n+1} - \delta_n)f(x_n) \\
&\quad + \delta_{n+1}(x_{n+2} - x_{n+1}) + (\delta_{n+1} - \delta_n)x_{n+1}\| \\
&= \|(1-\delta_{n+1})(f(x_{n+1}) - f(x_n)) + (\delta_{n+1} - \delta_n)(x_{n+1} - f(x_n)) \\
&\quad + \delta_{n+1}(x_{n+2} - x_{n+1})\| \leq (1-\delta_{n+1})\|f(x_{n+1}) - f(x_n)\| \\
&\quad + (\delta_{n+1} - \delta_n)\|x_{n+1} - f(x_n)\| + \delta_{n+1}\|x_{n+2} - x_{n+1}\| \\
&\leq (1-\delta_{n+1})\psi\|x_{n+1} - x_n\| + (\delta_{n+1} - \delta_n)\|x_{n+1} - f(x_n)\| \\
&\quad + \delta_{n+1}\|x_{n+2} - x_{n+1}\|.
\end{aligned}$$
(26)

**Theorem 15.** *Let K be a nonempty closed convex subset of a uniformly smooth Banach space E and $f : K \longrightarrow K$ a generalized contraction mapping. Let T be a nonexpansive self-mapping defined on K with $F(T) \neq \emptyset$. Assume that the conditions (i)-(v) of Assumption 12 are satisfied. Then, the iterative sequence $\{x_n\}_{n=1}^{\infty}$ which is defined from an arbitrary $x_1 \in K$ by (7) converges strongly to a fixed point p of T, which solves the variational inequality (6), given by*

$$\langle (I - f)p, J(x - p) \rangle \geq 0, \quad \text{for all } x \in F(T). \quad (27)$$

*Proof.* Setting $z_n = (x_{n+1} - \alpha_n^2 x_n)/(1 - \alpha_n^2)$ and $y_n = (1 - \delta_n)f(x_n) + \delta_n x_{n+1}$, one can obtain that

$$\begin{aligned}
&z_{n+1} - z_n \\
&= \frac{x_{n+2} - \alpha_{n+1}^2 x_{n+1}}{1 - \alpha_{n+1}^2} - \frac{x_{n+1} - \alpha_n^2 x_n}{1 - \alpha_n^2} \\
&= \frac{\alpha_{n+1}^1 f(x_{n+1}) + \alpha_{n+1}^3 T(y_{n+1})}{1 - \alpha_{n+1}^2} - \frac{\alpha_n^1 f(x_n) + \alpha_n^3 T(y_n)}{1 - \alpha_n^2} \\
&= \frac{\alpha_{n+1}^1}{1 - \alpha_{n+1}^2}(f(x_{n+1}) - f(x_n)) + \left(\frac{\alpha_{n+1}^1}{1 - \alpha_{n+1}^2} - \frac{\alpha_n^1}{1 - \alpha_n^2}\right) f(x_n) \\
&\quad + \frac{\alpha_{n+1}^3}{1 - \alpha_{n+1}^2}(T(y_{n+1}) - T(y_n)) + \left(\frac{\alpha_{n+1}^3}{1 - \alpha_{n+1}^2} - \frac{\alpha_n^3}{1 - \alpha_n^2}\right) T(y_n) \\
&= \frac{\alpha_{n+1}^1}{1 - \alpha_{n+1}^2}(f(x_{n+1}) - f(x_n)) - \left(\frac{\alpha_{n+1}^3}{1 - \alpha_{n+1}^2} - \frac{\alpha_n^3}{1 - \alpha_n^2}\right) f(x_n) \\
&\quad + \frac{\alpha_{n+1}^3}{1 - \alpha_{n+1}^2}(T(y_{n+1}) - T(y_n)) + \left(\frac{\alpha_{n+1}^3}{1 - \alpha_{n+1}^2} - \frac{\alpha_n^3}{1 - \alpha_n^2}\right) T(y_n) \\
&= \frac{\alpha_{n+1}^1}{1 - \alpha_{n+1}^2}(f(x_{n+1}) - f(x_n)) + \left(\frac{\alpha_{n+1}^3}{1 - \alpha_{n+1}^2} - \frac{\alpha_n^3}{1 - \alpha_n^2}\right) \\
&\quad \cdot (T(y_n) - f(x_n)) + \frac{\alpha_{n+1}^3}{1 - \alpha_{n+1}^2}(T(y_{n+1}) - T(y_n)).
\end{aligned}$$
(28)

Let $M_1 = \sup_n\{\|T(y_n) - f(x_n)\|\}$, $M_2 = \sup_n\{\|x_{n+1} - f(x_n)\|\}$, and $M = \max\{M_1, M_2\}$. Then,

$$\begin{aligned}
&\|z_{n+1} - z_n\| \\
&\leq \frac{\alpha_{n+1}^1}{1 - \alpha_{n+1}^2}\|f(x_{n+1}) - f(x_n)\| + \left|\frac{\alpha_{n+1}^3}{1 - \alpha_{n+1}^2} - \frac{\alpha_n^3}{1 - \alpha_n^2}\right| \\
&\quad \cdot \|T(y_n) - f(x_n)\| + \frac{\alpha_{n+1}^3}{1 - \alpha_{n+1}^2}\|T(y_{n+1}) - T(y_n)\| \\
&\leq \frac{\alpha_{n+1}^1}{1 - \alpha_{n+1}^2}\psi\|x_{n+1} - x_n\| + \left|\frac{\alpha_{n+1}^3}{1 - \alpha_{n+1}^2} - \frac{\alpha_n^3}{1 - \alpha_n^2}\right| \\
&\quad \cdot \|T(y_n) - f(x_n)\| + \frac{\alpha_{n+1}^3}{1 - \alpha_{n+1}^2}[(1 - \delta_{n+1})\psi\|x_{n+1} - x_n\| \\
&\quad + (\delta_{n+1} - \delta_n)\|x_{n+1} - f(x_n)\| + \delta_{n+1}\|x_{n+2} - x_{n+1}\|] \\
&\cdot (\text{by } (25)) = \frac{\alpha_{n+1}^1 \psi + \alpha_{n+1}^3(1 - \delta_{n+1})\psi}{1 - \alpha_{n+1}^2}\|x_{n+1} - x_n\| \\
&\quad + \left(\left|\frac{\alpha_{n+1}^3}{1 - \alpha_{n+1}^2} - \frac{\alpha_n^3}{1 - \alpha_n^2}\right| + \frac{\alpha_{n+1}^3(\delta_{n+1} - \delta_n)}{1 - \alpha_{n+1}^2}\right) M \\
&\quad + \frac{\alpha_{n+1}^3 \delta_{n+1}}{1 - \alpha_{n+1}^2}\|x_{n+2} - x_{n+1}\|.
\end{aligned}$$
(29)

It is needed to evaluate $\|x_{n+2} - x_{n+1}\|$.

$$\begin{aligned}
&x_{n+2} - x_{n+1} \\
&= \alpha_{n+1}^1 f(x_{n+1}) + \alpha_{n+1}^2 x_{n+1} + \alpha_{n+1}^3 Ty_{n+1} - \left(\alpha_n^1 f(x_n) \right. \\
&\quad \left. + \alpha_n^2 x_n + \alpha_n^3 Ty_n\right) = \alpha_{n+1}^1(f(x_{n+1}) - f(x_n)) \\
&\quad + \alpha_{n+1}^2(x_{n+1} - x_n) + \alpha_{n+1}^3(Ty_{n+1} - Ty_n) \\
&\quad + \left(\alpha_{n+1}^1 - \alpha_n^1\right)f(x_n) + \left(\alpha_{n+1}^2 - \alpha_n^2\right)x_n + \left(\alpha_{n+1}^3 - \alpha_n^3\right)Ty_n \\
&= \alpha_{n+1}^1(f(x_{n+1}) - f(x_n)) + \alpha_{n+1}^2(x_{n+1} - x_n) + \alpha_{n+1}^3 \\
&\quad \cdot (Ty_{n+1} - Ty_n) + \left(\left(\alpha_n^2 - \alpha_{n+1}^2\right) + \left(\alpha_n^3 - \alpha_{n+1}^3\right)\right)f(x_n) \\
&\quad + \left(\alpha_{n+1}^2 - \alpha_n^2\right)x_n + \left(\alpha_{n+1}^3 - \alpha_n^3\right)Ty_n = \alpha_{n+1}^1\left(f(x_{n+1}) \right. \\
&\quad \left. - f(x_n)\right) + \alpha_{n+1}^2(x_{n+1} - x_n) + \alpha_{n+1}^3(Ty_{n+1} - Ty_n) \\
&\quad + \left(\alpha_{n+1}^2 - \alpha_n^2\right)(x_n - f(x_n)) + \left(\alpha_{n+1}^3 - \alpha_n^3\right)(Ty_n - f(x_n)).
\end{aligned}$$
(30)

This leads to

$$\begin{aligned}
&\|x_{n+2} - x_{n+1}\| \\
&\leq \alpha_{n+1}^1 \psi\|x_{n+1} - x_n\| + \alpha_{n+1}^2\|x_{n+1} - x_n\| + \alpha_{n+1}^3\|Ty_{n+1} - Ty_n\| \\
&\quad + \left|\alpha_{n+1}^2 - \alpha_n^2\right|\|x_n - f(x_n)\| + \left|\alpha_{n+1}^3 - \alpha_n^3\right|\|Ty_n - f(x_n)\| \\
&\leq \alpha_{n+1}^1 \psi\|x_{n+1} - x_n\| + \alpha_{n+1}^2\|x_{n+1} - x_n\| + \alpha_{n+1}^3[(1 - \delta_{n+1})\psi \\
&\quad \cdot \|x_{n+1} - x_n\| + (\delta_{n+1} - \delta_n)\|x_{n+1} - f(x_n)\| \\
&\quad + \delta_{n+1}\|x_{n+2} - x_{n+1}\|] \text{ (by (25))} + \left|\alpha_{n+1}^2 - \alpha_n^2\right| \\
&\quad \cdot \|x_n - f(x_n)\| + \left|\alpha_{n+1}^3 - \alpha_n^3\right|\|Ty_n - f(x_n)\|
\end{aligned}$$

$$\begin{aligned}
&= \left(\alpha_{n+1}^2 + \left(\alpha_{n+1}^3 + \alpha_{n+1}^1\right)\psi - \alpha_{n+1}^3 \delta_{n+1}\psi\right)\|x_{n+1} - x_n\| \\
&\quad + \alpha_{n+1}^3 \delta_{n+1}\|x_{n+2} - x_{n+1}\| + \left(|\alpha_{n+1}^2 - \alpha_n^2| + |\alpha_{n+1}^3 - \alpha_n^3|\right. \\
&\quad \left. + \alpha_{n+1}^3 (\delta_{n+1} - \delta_n)\right)M = \left(\alpha_{n+1}^2 + (1 - \alpha_{n+1}^2)\psi - \alpha_{n+1}^3 \delta_{n+1}\psi\right) \\
&\quad \cdot \|x_{n+1} - x_n\| + \alpha_{n+1}^3 \delta_{n+1}\|x_{n+2} - x_{n+1}\| + \left(|\alpha_{n+1}^2 - \alpha_n^2|\right. \\
&\quad \left. + |\alpha_{n+1}^3 - \alpha_n^3| + \alpha_{n+1}^3 (\delta_{n+1} - \delta_n)\right)M = \left(\psi + \alpha_{n+1}^2 (1 - \psi)\right. \\
&\quad \left. - \alpha_{n+1}^3 \delta_{n+1}\psi\right)\|x_{n+1} - x_n\| + \alpha_{n+1}^3 \delta_{n+1}\|x_{n+2} - x_{n+1}\| \\
&\quad + \left(|\alpha_{n+1}^2 - \alpha_n^2| + |\alpha_{n+1}^3 - \alpha_n^3| + \alpha_{n+1}^3(\delta_{n+1} - \delta_n)\right)M \\
&= \left(\alpha_{n+1}^2 (1 - \psi) + (1 - \alpha_{n+1}^3 \delta_{n+1})\psi\right)\|x_{n+1} - x_n\| \\
&\quad + \alpha_{n+1}^3 \delta_{n+1}\|x_{n+2} - x_{n+1}\| + \left(|\alpha_{n+1}^2 - \alpha_n^2| + |\alpha_{n+1}^3 - \alpha_n^3|\right. \\
&\quad \left. + \alpha_{n+1}^3 (\delta_{n+1} - \delta_n)\right)M.
\end{aligned} \tag{31}$$

Let $d_n = (|\alpha_{n+1}^2 - \alpha_n^2| + |\alpha_{n+1}^3 - \alpha_n^3| + \alpha_{n+1}^3(\delta_{n+1} - \delta_n))$.
Therefore,

$$\begin{aligned}
\|x_{n+2} - x_{n+1}\| &\leq \frac{\alpha_{n+1}^2(1 - \psi) + (1 - \alpha_{n+1}^3 \delta_{n+1})\psi}{1 - \alpha_{n+1}^3 \delta_{n+1}}\|x_{n+1} - x_n\| \\
&\quad + \frac{d_n M}{1 - \alpha_{n+1}^3 \delta_{n+1}}.
\end{aligned} \tag{32}$$

Let $S_n = |(\alpha_{n+1}^3/(1 - \alpha_{n+1}^2)) - (\alpha_n^3/(1 - \alpha_n^2))| + ((\alpha_{n+1}^3(\delta_{n+1} - \delta_n))/(1 - \alpha_{n+1}^2))$ and substitute (32) into (29) to obtain

$$\begin{aligned}
\|z_{n+1} - z_n\| &\leq \left[\frac{\alpha_{n+1}^1 \psi + \alpha_{n+1}^3 (1 - \delta_{n+1})\psi}{1 - \alpha_{n+1}^2} + \frac{\alpha_{n+1}^3 \delta_{n+1}}{1 - \alpha_{n+1}^2} \times \frac{\alpha_{n+1}^2(1 - \psi) + (1 - \alpha_{n+1}^3 \delta_{n+1})\psi}{1 - \alpha_{n+1}^3 \delta_{n+1}}\right] \\
&\quad \cdot \|x_{n+1} - x_n\| + S_n M + \frac{\alpha_{n+1}^3 \delta_{n+1}}{1 - \alpha_{n+1}^2} \times \frac{d_n M}{1 - \alpha_{n+1}^3 \delta_{n+1}} \\
&= \left[\frac{\alpha_{n+1}^1 \psi + \alpha_{n+1}^3 (1 - \delta_{n+1})\psi - \alpha_{n+1}^3 \delta_{n+1}\left(\alpha_{n+1}^1 \psi + \alpha_{n+1}^3 (1 - \delta_{n+1})\psi\right)}{[1 - \alpha_{n+1}^2][1 - \alpha_{n+1}^3 \delta_{n+1}]} + \frac{\alpha_{n+1}^3 \delta_{n+1}\left(\alpha_{n+1}^2(1 - \psi) + (1 - \alpha_{n+1}^3 \delta_{n+1})\psi\right)}{[1 - \alpha_{n+1}^2][1 - \alpha_{n+1}^3 \delta_{n+1}]}\right] \\
&\quad \cdot \|x_{n+1} - x_n\| + \left(S_n + \frac{d_n \alpha_{n+1}^3 \delta_{n+1}}{[1 - \alpha_{n+1}^2][1 - \alpha_{n+1}^3 \delta_{n+1}]}\right)M \\
&= \left[\frac{\alpha_{n+1}^1 \psi + \alpha_{n+1}^3(1 - \delta_{n+1})\psi - \alpha_{n+1}^3 \delta_{n+1}\left(\alpha_{n+1}^1 \psi + \alpha_{n+1}^3 \psi - \alpha_{n+1}^3 \delta_{n+1}\psi\right)}{[1 - \alpha_{n+1}^2][1 - \alpha_{n+1}^3 \delta_{n+1}]} + \frac{\alpha_{n+1}^3 \delta_{n+1}\left(\alpha_{n+1}^2 - \alpha_{n+1}^2\psi + \psi - \alpha_{n+1}^3 \delta_{n+1}\psi\right)}{[1 - \alpha_{n+1}^2][1 - \alpha_{n+1}^3 \delta_{n+1}]}\right] \\
&\quad \cdot \|x_{n+1} - x_n\| + \left(S_n + \frac{d_n \alpha_{n+1}^3 \delta_{n+1}}{[1 - \alpha_{n+1}^2][1 - \alpha_{n+1}^3 \delta_{n+1}]}\right)M \\
&= \left[\frac{\alpha_{n+1}^1 \psi + \alpha_{n+1}^3(1 - \delta_{n+1})\psi - \alpha_{n+1}^3 \delta_{n+1}\left((1 - \alpha_{n+1}^2)\psi - \alpha_{n+1}^3 \delta_{n+1}\psi\right)}{[1 - \alpha_{n+1}^2][1 - \alpha_{n+1}^3 \delta_{n+1}]} + \frac{\alpha_{n+1}^3 \delta_{n+1}\left(\alpha_{n+1}^2 + (1 - \alpha_{n+1}^2)\psi - \alpha_{n+1}^3 \delta_{n+1}\psi\right)}{[1 - \alpha_{n+1}^2][1 - \alpha_{n+1}^3 \delta_{n+1}]}\right] \\
&\quad \cdot \|x_{n+1} - x_n\| + \left(S_n + \frac{d_n \alpha_{n+1}^3 \delta_{n+1}}{[1 - \alpha_{n+1}^2][1 - \alpha_{n+1}^3 \delta_{n+1}]}\right)M = \frac{\alpha_{n+1}^1 \psi + \alpha_{n+1}^3(1 - \delta_{n+1})\psi + \alpha_{n+1}^3 \delta_{n+1}\alpha_{n+1}^2}{[1 - \alpha_{n+1}^2][1 - \alpha_{n+1}^3 \delta_{n+1}]}\|x_{n+1} - x_n\| \\
&\quad + \left(S_n + \frac{d_n \alpha_{n+1}^3 \delta_{n+1}}{[1 - \alpha_{n+1}^2][1 - \alpha_{n+1}^3 \delta_{n+1}]}\right)M = \frac{(1 - \alpha_{n+1}^2)\psi - \alpha_{n+1}^3 \delta_{n+1}\psi + \alpha_{n+1}^3 \delta_{n+1}\alpha_{n+1}^2}{[1 - \alpha_{n+1}^2][1 - \alpha_{n+1}^3 \delta_{n+1}]}\|x_{n+1} - x_n\| \\
&\quad + \left(S_n + \frac{d_n \alpha_{n+1}^3 \delta_{n+1}}{[1 - \alpha_{n+1}^2][1 - \alpha_{n+1}^3 \delta_{n+1}]}\right)M = \left(1 - \frac{(1 - \alpha_{n+1}^2)(1 - \psi) - \alpha_{n+1}^3 \delta_{n+1}(1 - \psi)}{[1 - \alpha_{n+1}^2][1 - \alpha_{n+1}^3 \delta_{n+1}]}\right)\|x_{n+1} - x_n\| \\
&\quad + \left(S_n + \frac{d_n \alpha_{n+1}^3 \delta_{n+1}}{[1 - \alpha_{n+1}^2][1 - \alpha_{n+1}^3 \delta_{n+1}]}\right)M = \left(1 - \frac{(1 - \alpha_{n+1}^2)\phi - \alpha_{n+1}^3 \delta_{n+1}\phi}{[1 - \alpha_{n+1}^2][1 - \alpha_{n+1}^3 \delta_{n+1}]}\right)\|x_{n+1} - x_n\| \\
&\quad + \left(S_n + \frac{d_n \alpha_{n+1}^3 \delta_{n+1}}{[1 - \alpha_{n+1}^2][1 - \alpha_{n+1}^3 \delta_{n+1}]}\right)M = \left(1 - \frac{(1 - \alpha_{n+1}^2 - \alpha_{n+1}^3 \delta_{n+1})\phi}{[1 - \alpha_{n+1}^2][1 - \alpha_{n+1}^3 \delta_{n+1}]}\right)\|x_{n+1} - x_n\| \\
&\quad + \left(S_n + \frac{d_n \alpha_{n+1}^3 \delta_{n+1}}{[1 - \alpha_{n+1}^2][1 - \alpha_{n+1}^3 \delta_{n+1}]}\right)M \leq \left(1 - \frac{(1 - \alpha_{n+1}^2 - \alpha_{n+1}^3 \delta_{n+1})\phi}{1 - \alpha_{n+1}^2}\right)\|x_{n+1} - x_n\| \\
&\quad + \left(S_n + \frac{d_n \alpha_{n+1}^3 \delta_{n+1}}{[1 - \alpha_{n+1}^2][1 - \alpha_{n+1}^3 \delta_{n+1}]}\right)M.
\end{aligned} \tag{33}$$

It then follows that

$$\begin{aligned}&\|z_{n+1} - z_n\| - \|x_{n+1} - x_n\| \\ &\leq -\frac{(1 - \alpha_{n+1}^2 - \alpha_{n+1}^3 \delta_{n+1})\phi}{1 - \alpha_{n+1}^2} \|x_{n+1} - x_n\| \\ &\quad + \left(S_n + \frac{d_n \alpha_{n+1}^3 \delta_{n+1}}{(1 - \alpha_{n+1}^2)(1 - \alpha_{n+1}^3 \delta_{n+1})}\right) M,\end{aligned} \quad (34)$$

and thus,

$$\limsup_{n \to \infty} (\|z_{n+1} - z_n\| - \|x_{n+1} - x_n\|) \leq 0. \quad (35)$$

Invoking Lemma 8 to obtain that

$$\lim_{n \to \infty} \|z_n - x_n\| = 0. \quad (36)$$

Consequently,

$$\begin{aligned}\|x_{n+1} - x_n\| &= \|(1 - \alpha_n^2) z_n + \alpha_n^2 x_n - x_n\| \\ &= \|(1 - \alpha_n^2) z_n - (1 - \alpha_n^2) x_n\| \\ &= \|(1 - \alpha_n^2)(z_n - x_n)\| \\ &\leq (1 - \alpha_n^2) \|z_n - x_n\| \to 0 \text{ as } n \to \infty.\end{aligned} \quad (37)$$

Next is to show that $\lim_{n \to \infty} \|x_n - T(x_n)\| = 0$. From (7), it is obtained that

$$\begin{aligned}\|x_n - Tx_n\| &\leq \|x_n - x_{n+1}\| + \|x_{n+1} - T(x_n)\| \\ &\leq \|x_{n+1} - x_n\| + \|\alpha_n^1 f(x_n) + \alpha_n^2 x_n + \alpha_n^3 T(y_n) \\ &\quad - T(x_n)\| \leq \|x_{n+1} - x_n\| + \alpha_n^1 \|f(x_n) - T(x_n)\| \\ &\quad + \alpha_n^2 \|x_n - T(x_n)\| + \alpha_n^3 \|T(y_n) - T(x_n)\| \\ &\leq \|x_{n+1} - x_n\| + \alpha_n^1 \|f(x_n) - T(x_n)\| \\ &\quad + \alpha_n^2 \|x_n - T(x_n)\| + \alpha_n^3 \|y_n - x_n\| \\ &\leq \|x_{n+1} - x_n\| + \alpha_n^1 \|f(x_n) - T(x_n)\| \\ &\quad + \alpha_n^2 \|x_n - T(x_n)\| + \alpha_n^3 \|(1 - \delta_n) f(x_n) \\ &\quad + \delta_n x_{n+1} - x_n\| \leq \|x_{n+1} - x_n\| \\ &\quad + \alpha_n^1 \|f(x_n) - T(x_n)\| + \alpha_n^2 \|x_n - T(x_n)\| \\ &\quad + \alpha_n^3 (1 - \delta_n) \|x_n - f(x_n)\| + \alpha_n^3 \delta_n \|x_{n+1} - x_n\| \\ &= (1 + \alpha_n^3 \delta_n) \|x_{n+1} - x_n\| + (\alpha_n^1 + \alpha_n^3 (1 - \delta_n)) Q \\ &\quad + \alpha_n^2 \|x_n - T(x_n)\| = (1 + \alpha_n^3 \delta_n) \|x_{n+1} - x_n\| \\ &\quad + (1 - \alpha_n^3 \delta_n - \alpha_n^2) Q + \alpha_n^2 \|x_n - T(x_n)\|.\end{aligned} \quad (38)$$

Since $0 < \liminf_{n \to \infty} \alpha_n^2 \leq \limsup_{n \to \infty} \alpha_n^2 < 1$, let $0 < \eta \leq \alpha_n^2 < 1$, then

$$\begin{aligned}\|x_n - Tx_n\| &\leq \frac{1 + \alpha_n^3 \delta_n}{1 - \alpha_n^2} \|x_{n+1} - x_n\| + \frac{1 - \alpha_n^3 \delta_n - \alpha_n^2}{1 - \alpha_n^2} Q \\ &\leq \frac{1 + \alpha_n^3 \delta_n}{1 - \eta} \|x_{n+1} - x_n\| + \frac{1 - \alpha_n^3 \delta_n - \alpha_n^2}{1 - \eta} Q,\end{aligned} \quad (39)$$

which goes to zero as $n \to \infty$ by (37) and condition (ii) of Assumption 12.

Let a sequence $\{x_t\}$ be defined by $x_t = t f(x_t) + (1 - t) T x_t$ for $t \in (0, 1)$. It is known by Lemma 10 that $\{x_t\}$ converges strongly to $p \in F(T)$, which solves the variational inequality:

$$\langle f(p) - p, J(x - p) \rangle \leq 0, \quad \forall x \in F(T), \quad (40)$$

which is equivalent to

$$\langle (I - f) p, J(x - p) \rangle \geq 0, \quad \forall x \in F(T). \quad (41)$$

It can be shown that

$$\limsup_{n \to \infty} \langle f(p) - p, J(x_{n+1} - p) \rangle \leq 0, \quad (42)$$

where $p \in F(T)$ is the unique fixed point of the generalized contraction $P_{F(T)} f(p)$ (Proposition 7), that is, $p = P_{F(T)} f(p)$.

By (39), $\lim_{n \to \infty} \|x_n - Tx_n\| = 0$, it follows from Lemma 9 that

$$\limsup_{n \to \infty} \langle f(p) - p, J(x_n - p) \rangle \leq 0. \quad (43)$$

Due to the continuity of the duality map and the fact that $\|x_{n+1} - x_n\| \to 0$ as $n \to \infty$ by (37), it is obtained that

$$\begin{aligned}&\limsup_{n \to \infty} \langle f(p) - p, J(x_{n+1} - p) \rangle \\ &= \limsup_{n \to \infty} \langle f(p) - p, J(x_{n+1} - x_n + x_n - p) \rangle \\ &= \limsup_{n \to \infty} \langle f(p) - p, J(x_n - p) \rangle \leq 0.\end{aligned} \quad (44)$$

Lastly, it is shown that $x_n \to p \in F(T)$ as $n \to \infty$.

Suppose that the sequence $\{x_n\}_{n=1}^{\infty}$ does not converge strongly to $p \in F(T)$. Then, there exists $\varepsilon > 0$ and a subsequence $\{x_{n_k}\}_{k=1}^{\infty}$ of $\{x_n\}_{n=1}^{\infty}$ such that $\|x_{n_k} - p\| \geq \varepsilon$, for all $k \in \mathbb{N}$. Therefore, for this $\varepsilon$, there exists $c \in (0, 1/2)$ such that

$$\|f(x_{n_k}) - f(p)\| \leq c \|x_{n_k} - p\|. \quad (45)$$

$$\begin{aligned}
&\|x_{n_{k+1}} - p\|^2 \\
&= \alpha_{n_k}^1 \langle f(x_{n_k}) - p, J(x_{n_{k+1}} - p) \rangle + \alpha_{n_k}^2 \langle x_{n_k} - p, J(x_{n_{k+1}} - p) \rangle \\
&\quad + \alpha_{n_k}^3 \langle T(y_{n_k}) - p, J(x_{n_{k+1}} - p) \rangle = \alpha_{n_k}^1 \langle f(x_{n_k}) \\
&\quad - f(p), J(x_{n_{k+1}} - p) \rangle + \alpha_n^1 \langle f(p) - p, J(x_{n_{k+1}} - p) \rangle + \alpha_{n_k}^2 \\
&\quad \cdot \langle x_{n_k} - p, J(x_{n_{k+1}} - p) \rangle + \alpha_{n_k}^3 \langle T(y_{n_k}) - p, J(x_{n_{k+1}} - p) \rangle \\
&\leq c\alpha_{n_k}^1 \|x_{n_k} - p\| \|x_{n_{k+1}} - p\| + \alpha_n^1 \langle f(p) - p, J(x_{n_{k+1}} - p) \rangle \\
&\quad + \alpha_{n_k}^2 \|x_{n_k} - p\| \|x_{n_{k+1}} - p\| + \alpha_{n_k}^3 \|(1 - \delta_{n_k}) f(x_{n_k}) \\
&\quad + \delta_{n_k} x_{n_{k+1}} - p\| \|x_{n_{k+1}} - p\| \leq c\alpha_{n_k}^1 \|x_{n_k} - p\| \|x_{n_{k+1}} - p\| \\
&\quad + \alpha_n^1 \langle f(p) - p, J(x_{n_{k+1}} - p) \rangle + \alpha_{n_k}^2 \|x_{n_k} - p\| \|x_{n_{k+1}} - p\| \\
&\quad + \alpha_{n_k}^3 (1 - \delta_{n_k}) \|f(x_{n_k}) - p\| \|x_{n_{k+1}} - p\| \\
&\quad + \alpha_{n_k}^3 \delta_{n_k} \|x_{n_{k+1}} - p\|^2 \leq c\alpha_{n_k}^1 \|x_{n_k} - p\| \|x_{n_{k+1}} - p\| \\
&\quad + \alpha_n^1 \langle f(p) - p, J(x_{n_{k+1}} - p) \rangle + \alpha_{n_k}^2 \|x_{n_k} - p\| \|x_{n_{k+1}} - p\| \\
&\quad + c\alpha_{n_k}^3 (1 - \delta_{n_k}) \|x_{n_k} - p\| \|x_{n_{k+1}} - p\| + \alpha_{n_k}^3 (1 - \delta_{n_k}) \\
&\quad \cdot \|f(p) - p\| \|x_{n_{k+1}} - p\| + \alpha_{n_k}^3 \delta_{n_k} \|x_{n_{k+1}} - p\|^2 \\
&= \left( c\alpha_{n_k}^1 + \alpha_{n_k}^2 + c\alpha_{n_k}^3 (1 - \delta_{n_k}) \right) \|x_{n_k} - p\| \|x_{n_{k+1}} - p\| \\
&\quad + \alpha_n^1 \langle f(p) - p, J(x_{n_{k+1}} - p) \rangle + \alpha_{n_k}^3 (1 - \delta_{n_k}) \|f(p) - p\| \\
&\quad \cdot \|x_{n_{k+1}} - p\| + \alpha_{n_k}^3 \delta_{n_k} \|x_{n_{k+1}} - p\|^2 \leq \frac{1}{2} \bigl( c\alpha_{n_k}^1 + \alpha_{n_k}^2 \\
&\quad + c\alpha_{n_k}^3 (1 - \delta_{n_k}) \bigr) \left( \|x_{n_k} - p\|^2 + \|x_{n_{k+1}} - p\|^2 \right) + \alpha_n^1 \langle f(p) \\
&\quad - p, J(x_{n_{k+1}} - p) \rangle + \alpha_{n_k}^3 \delta_{n_k} \|x_{n_{k+1}} - p\|^2 + \frac{1}{2} \alpha_{n_k}^3 (1 - \delta_{n_k}) \\
&\quad \cdot \left( \|f(p) - p\|^2 + \|x_{n_{k+1}} - p\|^2 \right) = \frac{1}{2} \bigl( c (\alpha_{n_k}^1 + \alpha_{n_k}^3 \\
&\quad \cdot (1 - \delta_{n_k})) + \alpha_{n_k}^2 \bigr) \|x_{n_k} - p\|^2 + \alpha_n^1 \langle f(p) - p, J(x_{n_{k+1}} - p) \rangle \\
&\quad + \frac{1}{2} \bigl( c(\alpha_{n_k}^1 + \alpha_{n_k}^3 (1 - \delta_{n_k})) + \alpha_{n_k}^2 + 2\alpha_{n_k}^3 \delta_{n_k} \\
&\quad + \alpha_{n_k}^3 (1 - \delta_{n_k}) \bigr) \|x_{n_{k+1}} - p\|^2 + \frac{1}{2} \alpha_{n_k}^3 (1 - \delta_{n_k}) \\
&\quad \cdot \|f(p) - p\|^2 = \frac{1}{2} \bigl( c(\alpha_{n_k}^1 + \alpha_{n_k}^3 (1 - \delta_{n_k})) + \alpha_{n_k}^2 \bigr) \\
&\quad \cdot \|x_{n_k} - p\|^2 + \alpha_n^1 \langle f(p) - p, J(x_{n_{k+1}} - p) \rangle \\
&\quad + \frac{1}{2} \bigl( c(\alpha_{n_k}^1 + \alpha_{n_k}^3 (1 - \delta_{n_k})) + \alpha_{n_k}^2 + \alpha_{n_k}^3 (1 + \delta_{n_k}) \bigr) \\
&\quad \cdot \|x_{n_{k+1}} - p\|^2 + \frac{1}{2} \alpha_{n_k}^3 (1 - \delta_{n_k}) \|f(p) - p\|^2 \\
&= \frac{1}{2} \bigl( c(1 - \alpha_{n_k}^2 - \alpha_{n_k}^3 \delta_{n_k}) + \alpha_{n_k}^2 \bigr) \|x_{n_k} - p\|^2 \\
&\quad + \alpha_n^1 \langle f(p) - p, J(x_{n_{k+1}} - p) \rangle + \frac{1}{2} \bigl( c(1 - \alpha_{n_k}^2 - \alpha_{n_k}^3 \delta_{n_k}) \\
&\quad + \alpha_{n_k}^2 + \alpha_{n_k}^3 (1 + \delta_{n_k}) \bigr) \|x_{n_{k+1}} - p\|^2 + \frac{1}{2} \alpha_{n_k}^3 (1 - \delta_{n_k}) \\
&\quad \cdot \|f(p) - p\|^2.
\end{aligned}$$
(46)

Observe that

$$\begin{aligned}
&2 - c(1 - \alpha_{n_k}^2 - \alpha_{n_k}^3 \delta_{n_k}) - \alpha_{n_k}^2 - \alpha_{n_k}^3 (1 + \delta_{n_k}) \\
&= 2 - c + c\alpha_{n_k}^2 + c\alpha_{n_k}^3 \delta_{n_k} - \alpha_{n_k}^2 - \alpha_{n_k}^3 - \alpha_{n_k}^3 \delta_{n_k} \\
&= 2 - c - (1-c)\alpha_{n_k}^2 - (1-c)\alpha_{n_k}^3 \delta_{n_k} - \alpha_{n_k}^3 \\
&= 1 - c - (1-c)\alpha_{n_k}^2 - (1-c)\alpha_{n_k}^3 \delta_{n_k} + 1 - \alpha_{n_k}^3 \\
&= 1 + (1-c)\left(1 - \alpha_{n_k}^2 - \alpha_{n_k}^3 \delta_{n_k}\right) - \alpha_{n_k}^3,
\end{aligned}$$
(47)

$$\alpha_{n_k}^1 = 1 - \alpha_{n_k}^2 - \alpha_{n_k}^3 \leq 1 - \alpha_{n_k}^2 - \alpha_{n_k}^3 \delta_{n_k} \text{ (since } \delta_{n_k} \in (0,1)\text{)}.$$
(48)

Multiplying (46) by 2 gives

$$\begin{aligned}
&\|x_{n_{k+1}} - p\|^2 \\
&\leq \frac{c(1 - \alpha_{n_k}^2 - \alpha_{n_k}^3 \delta_{n_k}) + \alpha_{n_k}^2}{1 + (1-c)(1 - \alpha_{n_k}^2 - \alpha_{n_k}^3 \delta_{n_k}) - \alpha_{n_k}^3} \|x_{n_k} - p\|^2 \\
&\quad + \frac{\alpha_n^1}{1 + (1-c)(1 - \alpha_{n_k}^2 - \alpha_{n_k}^3 \delta_{n_k}) - \alpha_{n_k}^3} \langle f(p) - p, J(x_{n_{k+1}} - p) \rangle \\
&\quad + \frac{\alpha_{n_k}^3 (1 - \delta_{n_k})}{1 + (1-c)(1 - \alpha_{n_k}^2 - \alpha_{n_k}^3 \delta_{n_k}) - \alpha_{n_k}^3} \|f(p) - p\|^2 \\
&= \left( 1 - \frac{(1-2c)(1 - \alpha_{n_k}^2 - \alpha_{n_k}^3 \delta_{n_k}) + \alpha_{n_k}^1}{1 + (1-c)(1 - \alpha_{n_k}^2 - \alpha_{n_k}^3 \delta_{n_k}) - \alpha_{n_k}^3} \right) \|x_{n_k} - p\|^2 \\
&\quad + \frac{\alpha_{n_k}^1}{1 + (1-c)(1 - \alpha_{n_k}^2 - \alpha_{n_k}^3 \delta_{n_k}) - \alpha_{n_k}^3} \langle f(p) - p, J(x_{n_{k+1}} - p) \rangle \\
&\quad + \frac{\alpha_{n_k}^3 (1 - \delta_{n_k})}{1 + (1-c)(1 - \alpha_{n_k}^2 - \alpha_{n_k}^3 \delta_{n_k}) - \alpha_{n_k}^3} \|f(p) - p\|^2 \\
&\leq \left( 1 - \frac{(1-2c)(1 - \alpha_{n_k}^2 - \alpha_{n_k}^3 \delta_{n_k})}{1 + (1-c)(1 - \alpha_{n_k}^2 - \alpha_{n_k}^3 \delta_{n_k}) - \alpha_{n_k}^3} \right) \|x_{n_k} - p\|^2 \\
&\quad + \frac{(1-2c)(1 - \alpha_{n_k}^2 - \alpha_{n_k}^3 \delta_{n_k})}{1 + (1-c)(1 - \alpha_{n_k}^2 - \alpha_{n_k}^3 \delta_{n_k}) - \alpha_{n_k}^3} \cdot \frac{1}{1 - 2c} \\
&\quad \cdot \langle f(p) - p, J(x_{n_{k+1}} - p) \rangle \\
&\quad + \frac{\alpha_{n_k}^3 (1 - \delta_{n_k})}{1 + (1-c)(1 - \alpha_{n_k}^2 - \alpha_{n_k}^3 \delta_{n_k}) - \alpha_{n_k}^3} \\
&\quad \cdot \|f(p) - p\|^2 \text{ (by (48))}.
\end{aligned}$$
(49)

Using Lemma 11, it shows that $x_{n_k} \longrightarrow p$ as $k \longrightarrow \infty$. A contradiction, hence, $\{x_n\}_{n=1}^{\infty}$ converges strongly to $p \in F(T)$.

The next result shows that under suitable conditions, the implicit iterative sequences (5) and (7) converge to the same fixed point.

**Theorem 16.** *Let K be a nonempty closed convex subset of a uniformly smooth Banach space E and $f : K \longrightarrow K$ a c-contraction mapping with $c \in [0, 1)$. Let T be a nonexpansive self-mapping defined on K with $F(T) \neq \emptyset$. Let $\{\{\alpha_n^i\}_{n=1}^{\infty}\}_{i=1}^{3} \subset [0, 1]$ and $\{\delta_n\}_{n=1}^{\infty} \subset (0, 1)$ be real sequences such that $\sum_{i=1}^{3} \alpha_n^i = 1$. Given that $\lim_{n \to \infty} \alpha_n^3/(1 - \alpha_n^2 - \alpha_n^3 \delta_n) = 0$, then $\{x_n\}_{n=1}^{\infty}$ defined by (7) converges to p if and only if $\{y_n\}_{n=1}^{\infty}$ defined by (5) converges to p.*

*Proof.* Notice that (7) and (5) are, respectively, given by

$$\begin{aligned}
x_{n+1} &= \alpha_n^1 f(x_n) + \alpha_n^2 x_n \\
&\quad + \alpha_n^3 T((1-\delta_n)f(x_n) + \delta_n x_{n+1}), \quad n \in \mathbb{N}, \\
y_{n+1} &= \alpha_n f(y_n) + \beta_n y_n \\
&\quad + \gamma_n T(\delta_n y_n + (1-\delta_n)y_{n+1}), \quad n \in \mathbb{N}.
\end{aligned} \quad (50)$$

It is needed to show that $\|x_n - y_n\| \longrightarrow 0$, as $n \longrightarrow \infty$.

$$\begin{aligned}
&\|x_{n+1} - y_{n+1}\| \\
&= \|\alpha_n^1 f(x_n) + \alpha_n^2 x_n + \alpha_n^3 T((1-\delta_n)f(x_n) + \delta_n x_{n+1}) \\
&\quad - (\alpha_n^1 f(y_n) + \alpha_n^2 y_n + \alpha_n^3 T(\delta_n y_n + (1-\delta_n)y_{n+1}))\| \\
&= \|\alpha_n^1(f(x_n) - f(y_n)) + \alpha_n^2(x_n - y_n) + \alpha_n^3(T((1-\delta_n)f(x_n) \\
&\quad + \delta_n x_{n+1}) - T(\delta_n y_n + (1-\delta_n)y_{n+1}))\| \\
&\leq \alpha_n^1 \|f(x_n) - f(y_n)\| + \alpha_n^2 \|x_n - y_n\| + \alpha_n^3 \|T((1-\delta_n)f(x_n) \\
&\quad + \delta_n x_{n+1}) - T(\delta_n y_n + (1-\delta_n)y_{n+1})\| \\
&\leq \alpha_n^1 c \|x_n - y_n\| + \alpha_n^2 \|x_n - y_n\| + \alpha_n^3 \|(1-\delta_n)(f(x_n) - y_{n+1}) \\
&\quad + \delta_n(x_{n+1} - y_n)\| \leq \alpha_n^1 c \|x_n - y_n\| + \alpha_n^2 \|x_n - y_n\| \\
&\quad + \alpha_n^3 (1-\delta_n) \|f(x_n) - f(y_n) + f(y_n) - y_{n+1}\| \\
&\quad + \alpha_n^3 \delta_n \|x_{n+1} - y_{n+1} + y_{n+1} - y_n\| \\
&\leq \alpha_n^1 c \|x_n - y_n\| + \alpha_n^2 \|x_n - y_n\| + \alpha_n^3 (1-\delta_n) c \|x_n - y_n\| \\
&\quad + \alpha_n^3 (1-\delta_n) \|y_{n+1} - f(y_n)\| + \alpha_n^3 \delta_n \|x_{n+1} - y_{n+1}\| \\
&\quad + \alpha_n^3 \delta_n \|y_{n+1} - y_n\| = (\alpha_n^1 c + \alpha_n^3 (1-\delta_n) c + \alpha_n^2) \|x_n - y_n\| \\
&\quad + \alpha_n^3 \delta_n \|x_{n+1} - y_{n+1}\| + \alpha_n^3 (1-\delta_n) \|y_{n+1} - f(y_n)\| \\
&\quad + \alpha_n^3 \delta_n \|y_{n+1} - y_n\|.
\end{aligned} \quad (51)$$

Since $\{y_n\}_{n=1}^{\infty}$ and $\{f(y_n)\}_{n=1}^{\infty}$ are bounded, let $M_2 = \max\{\sup_n \|y_{n+1} - f(y_n)\|, \sup_n \|y_{n+1} - y_n\|\}$. Then,

$$\begin{aligned}
\|x_{n+1} - y_{n+1}\| &\leq \frac{\alpha_n^1 c + \alpha_n^3 (1-\delta_n) c + \alpha_n^2}{1 - \alpha_n^3 \delta_n} \|x_n - y_n\| \\
&\quad + \frac{\alpha_n^3}{1 - \alpha_n^3 \delta_n} M_2 = (1 - \beta_n) \|x_n - y_n\| \\
&\quad + \frac{\alpha_n^3}{(1 - \alpha_n^2 - \alpha_n^3 \delta_n)(1-c)} \beta_n M_2,
\end{aligned} \quad (52)$$

where $\beta_n = ((1 - \alpha_n^2 - \alpha_n^3 \delta_n)(1-c))/(1 - \alpha_n^3 \delta_n)$. From the given condition, it follows that.

$\limsup_{n \to \infty} \alpha_n^3/(1 - \alpha_n^2 - \alpha_n^3 \delta_n) \leq 0$. Apply Lemma 11 to (52) and take $\gamma_n = 0$ to get that $\|x_n - y_n\| \longrightarrow 0$, as $n \longrightarrow \infty$. Next, suppose $\|y_n - p\| \longrightarrow 0$ as $n \longrightarrow \infty$. It follows that

$$\begin{aligned}
\|x_n - p\| &= \|x_n - y_n + y_n - p\| \\
&\leq \|x_n - y_n\| + \|y_n - p\| \longrightarrow 0 \text{ as } n \longrightarrow \infty.
\end{aligned} \quad (53)$$

Similarly, suppose $\|x_n - p\| \longrightarrow 0$ as $n \longrightarrow \infty$. Then,

$$\begin{aligned}
\|y_n - p\| &= \|y_n - x_n + x_n - p\| \\
&\leq \|y_n - x_n\| + \|x_n - p\| \longrightarrow 0 \text{ as } n \longrightarrow \infty.
\end{aligned} \quad (54)$$

Hence, the implicit iterative sequences (5) and (7) converge to the same fixed point under suitable conditions.

*Remark 17.* One can deduce the following results from Theorem 15.

**Corollary 18.** *Let K be a nonempty closed convex subset of a uniformly smooth Banach space E and T a nonexpansive self-mapping defined on K with $F(T) \neq \emptyset$. Assume that the real sequences $\{\alpha_n\}_{n=1}^{\infty} \subset (0, 1)$ and $\{\delta_n\}_{n=1}^{\infty} \subset (0, 1)$ satisfy the conditions:*

(i) $\lim_{n \to \infty} \alpha_n = 0$

(ii) $\sum_{n=1}^{\infty} \alpha_n = \infty$

(iii) $\sum_{n=1}^{\infty} |\alpha_{n+1} - \alpha_n| < \infty$

(iv) $0 < \varepsilon \leq \delta_n \leq \delta_{n+1} < 1$

*Then, the iterative sequence $\{x_n\}_{n=1}^{\infty}$ which is defined from an arbitrary $x_1 \in K$ by*

$$x_{n+1} = \alpha_n x_n + (1 - \alpha_n) T((1-\delta_n)x_n + \delta_n x_{n+1}) \quad (55)$$

*converges strongly to a fixed point p of T which solves the variational inequality (6).*

*Proof.* The result follows from Theorem 15 by simply taking $f$ to be the identity mapping on $K$.

**Corollary 19.** *Let K be a nonempty closed convex subset of a uniformly smooth Banach space E and T a nonexpansive self-mapping defined on K with $F(T) \neq \emptyset$. Assume that the real sequence $\{\alpha_n\}_{n=1}^{\infty} \subset (0, 1)$ satisfies the following conditions:*

(i) $\lim_{n \to \infty} \alpha_n = 0$

(ii) $\sum_{n=1}^{\infty} \alpha_n = \infty$

(iii) $\sum_{n=1}^{\infty} |\alpha_{n+1} - \alpha_n| < \infty$

*Then, the iterative sequence $\{x_n\}_{n=1}^{\infty}$ which is defined from an arbitrary $x_1 \in K$ by*

$$x_{n+1} = \alpha_n x_n + (1 - \alpha_n) T\left(\frac{x_n + x_{n+1}}{2}\right) \quad (56)$$

*converges strongly to a fixed point p of T which solves the variational inequality (6).*

*Proof.* The result follows from Theorem 15 by simply taking $f$ to be the identity mapping on $K$ and $\delta_n = 1/2$ for all $n \in \mathbb{N}$. Consequently, this improves and extend the results of Alghamdi et al. [19].

## 4. Applications

*4.1. Application to Fixed Points of $\lambda$-Strictly Pseudocontractive Mappings.* Let $K$ be a closed convex subset of a real Banach space $E$. A mapping $S : K \longrightarrow K$ is said to be $\lambda$-strictly pseudocontractive mapping if there exists $0 \leq \lambda < 1$ such that

$$\|Sx - Sy\|^2 \leq \|x - y\|^2 - \lambda\|(I - S)x - (I - S)y\|^2, \quad \forall x, y \in K, \tag{57}$$

where $I$ denotes the identity operator on $K$.

Zhou [20] established the following lemma which gives a relationship between $\lambda$-strictly pseudocontractive mappings and nonexpansive mappings.

**Lemma 20.** *Let $K$ be a nonempty subset of a 2-uniformly smooth Banach space $E$. Let $S : K \longrightarrow K$ be a $\lambda$-strictly pseudocontractive mapping. For $\theta \in (0, 1)$, define*

$$Tx = \theta x + (1 - \theta)Sx, \quad \forall x \in K. \tag{58}$$

*Then, as $\theta \in (0, \lambda/L^2]$ (where $L$ is the 2-uniformly smooth constant of a 2-uniformly smooth Banach space), $T : K \longrightarrow K$ is nonexpansive such that $F(T) = F(S)$.*

The following result is obtained by using Lemma 20 and Theorem 15.

**Corollary 21.** *Let $K$ be a nonempty closed convex subset of a 2-uniformly smooth Banach space $E$ and $f : K \longrightarrow K$ a generalized contraction mapping. Let $S : K \longrightarrow K$ a $\lambda$-strictly pseudocontractive mapping with $F(T) \neq \emptyset$. Suppose that the conditions (i)-(v) of Assumption 12 are satisfied and $T$ is a mapping from $K$ into itself, defined by $Tx = \alpha x + (1 - \theta)Sx$, $x \in K, \theta \in (0, 1)$. Then, for an arbitrary $x_1 \in K$, the iterative sequence $\{x_n\}_{n=1}^{\infty}$ defined by*

$$\begin{aligned} x_{n+1} &= \alpha_n^1 f(x_n) + \alpha_n^2 x_n \\ &\quad + \alpha_n^3 T((1 - \delta_n)f(x_n) + \delta_n x_{n+1}), \quad \text{for all } n \in \mathbb{N}, \end{aligned} \tag{59}$$

*converges strongly to a fixed point $p$ of $S$, which solves the variational inequality*

$$\langle (I - f)p, J(x - p) \rangle \geq 0, \quad \text{for all } x \in F(S). \tag{60}$$

*4.2. Application to Solution of $\alpha$-Inverse-Strongly Monotone Mappings.* Let $K$ be a nonempty closed convex subset of a Hilbert space $H$. The metric projection $P_K$ is defined from $H$ onto $K$ by

$$P_K x := \arg \min_{y \in K} \|x - y\|^2, \quad x \in H \tag{61}$$

and characterized by

$$P_K(x) := \arg \min_{z \in K} \|x - z\|^2, \quad x \in H. \tag{62}$$

$P_K(x)$ is known as the only point in $K$ that minimizes the objective $\|x - z\|$ over $z \in K$. A mapping $A$ of $K$ into $H$ is called monotone if $\langle Au - Av, u - v \rangle \geq 0$, for all $u, v \in K$. The classical variational inequality (VI) problem is to find $u^* \in K$ such that

$$\langle Au^*, u - u^* \rangle \geq 0, \quad u \in K, \tag{63}$$

where $A$ is a (single-valued) monotone operator in Hilbert space $H$ [21, 22]. In this work, the solution set of (63) is denoted by $\text{VI}(K, A)$. In the context of the variational inequality problem, (62) implies that

$$u \in \text{VI}(K, A) \Leftrightarrow u = P_K(u - \gamma Au), \quad \forall \gamma > 0. \tag{64}$$

$A$ is said to be $\alpha$-inverse-strongly monotone if there exists a positive real number $\alpha$ such that

$$\langle Au - Av, u - v \rangle \geq \alpha \|Au - Av\|^2, \tag{65}$$

for all $u, v \in K$. If $A$ is an $\alpha$-inverse-strongly monotone mapping of $K$ to $H$, it is known that $A$ is $1/\alpha$-Lipschitz continuous. Also, we have that for all $u, v \in K$ and $\gamma > 0$,

$$\begin{aligned} \|(I - \gamma A)u &- (I - \gamma A)v\|^2 \\ &= \|(u - v) - (Au - Av)\|^2 \\ &= \|u - v\|^2 - 2\gamma \langle u - v, Au - Av \rangle + \gamma^2 \|Au - Av\|^2 \\ &\leq \|u - v\|^2 + \gamma(\gamma - 2\alpha)\|Au - Av\|^2. \end{aligned} \tag{66}$$

Therefore, if $\gamma \leq 2\alpha$, then $I - \gamma A$ is a nonexpansive mapping of $K$ into $K$. Consequently, one can apply Theorem 15 to deduce the following result:

**Corollary 22.** *Let $K$ be a nonempty closed convex subset of a real Hilbert space $H$ and $f : K \longrightarrow K$ a generalized contractions. Let $A$ be an $\alpha$-inverse-strongly monotone mapping of $K$ to $H$ with $A^{-1}0 \neq \emptyset$. Assume that the conditions (i)-(v) of Assumption 12 are satisfied. Then, the iterative sequence $\{x_n\}_{n=1}^{\infty}$ which is defined from an arbitrary $x_1 \in K$ by*

$$\begin{aligned} x_{n+1} &= \alpha_n^1 f(x_n) + \alpha_n^2(x_n) + \alpha_n^3 P_K(I - \gamma A) \\ &\quad \cdot ((1 - \delta_n)f(x_n) + \delta_n x_{n+1}), \quad n \in \mathbb{N}, \end{aligned} \tag{67}$$

Journal of Applied Mathematics 11converges strongly to a solution $p$ in $A^{-1}0$, which solves the variational inequality

$$\langle (I-f)p, x-p \rangle \geq 0, \quad \text{for all } x \in A^{-1}0. \tag{68}$$

### 4.3. Application to Fredholm Integral Equation in Hilbert Spaces.
Consider a Fredholm integral equation of the form

$$x(t) = g(t) + \int_0^1 \Phi(t, s, x(s)) ds, \quad t \in [0,1], \tag{69}$$

where $g$ is a continuous function on $[0,1]$ and $\Phi : [0,1] \times [0,1] \times \mathbb{R} \longrightarrow \mathbb{R}$ is continuous. The existence of solutions of (69) has been studied (see [23] and the references therein). If $\Phi$ satisfies the Lipschitz continuity condition

$$|\Phi(t,s,x) - \Phi(t,s,y)| \leq |x - y|, \quad s,t \in [0,1], x,y \in \mathbb{R}, \tag{70}$$

then equation (69) has at least one solution in the Hilbert space $L^2[0,1]$ ([23], Theorem 3.3). Precisely, define a mapping $T : L^2[0,1] \longrightarrow L^2[0,1]$ by

$$Tx(t) = g(t) + \int_0^1 \Phi(t, s, x(s)) ds, \quad t \in [0,1]. \tag{71}$$

It is known that $T$ is nonexpansive. Indeed, for $x, y \in L^2[0,1]$

$$\begin{aligned} \|Tx - Ty\|^2 &= \int_0^1 |Tx(t) - Ty(t)|^2 dt \\ &= \int_0^1 \left| \int_0^1 \Phi(t,s,x(s)) - \Phi(t,s,y(s)) ds \right|^2 dt \\ &\leq \int_0^1 \left| \int_0^1 |x(s) - y(s)| ds \right|^2 dt \\ &\leq \int_0^1 |x(s) - y(s)|^2 ds = \|x - y\|^2. \end{aligned} \tag{72}$$

Thus, finding a solution of integral equation (69) is reduced to finding a fixed point of the nonexpansive mapping $T$ in the Hilbert space $L^2[0,1]$. Consequently, the following result is obtainable.

**Corollary 23.** *Let $K$ be a nonempty closed convex subset of a Hilbert space $L^2[0,1]$, $T : K \longrightarrow K$, defined by (71) with $F(T) \neq \emptyset$ and $f : K \longrightarrow K$ is a generalized contraction. Suppose that the conditions (i)-(v) of Assumption 12 are satisfied. Then, for an arbitrary $x_1 \in K$, the iterative sequence $\{x_n\}_{n=1}^{\infty}$ defined by*

$$\begin{aligned} x_{n+1} &= \alpha_n^1 f(x_n) + \alpha_n^2(x_n) \\ &\quad + \alpha_n^3 T((1 - \delta_n) f(x_n) + \delta_n x_{n+1}), \quad n \in \mathbb{N}, \end{aligned} \tag{73}$$

*converges strongly to a fixed point $p$ of $T$, which solves the variational inequality*

$$\langle (I-f)p, x-p \rangle \geq 0, \quad \text{for all } x \in F(T). \tag{74}$$

Examples of real sequences which satisfy the conditions of Assumption 12 are

$$\begin{aligned} \{\alpha_n^1\}_{n=1}^{\infty} &= \left\{ \frac{1}{2n} \right\}_{n=1}^{\infty}; \\ \{\alpha_n^2\}_{n=1}^{\infty} &= \left\{ 1 - \frac{3}{2n} \right\}_{n=1}^{\infty}; \\ \{\alpha_n^3\}_{n=1}^{\infty} &= \left\{ \frac{1}{n} \right\}_{n=1}^{\infty}; \\ \{\delta_n\}_{n=1}^{\infty} &= \left\{ \frac{n}{2(n+1)} \right\}_{n=1}^{\infty}. \end{aligned} \tag{75}$$

## Data Availability

Data sharing not applicable to this article as no datasets were generated or analyzed during the current study.

## Conflicts of Interest

The authors declare no conflict of interest.

## Authors' Contributions

All authors contributed significantly in writing this article. All authors read and approved the final manuscript.